\documentclass[twoside,reqno]{lt7-proc}
\usepackage{epsfig,cite}
\usepackage{amssymb,amsmath}
\usepackage{times}
 \newcommand{\abs}[1]{\lvert#1\rvert}
 \newcommand{\inv}{^{-1}}
\newcommand{\ten}{\otimes}

\begin{document}
\sloppy \raggedbottom
\setcounter{page}{1}

\newpage
\setcounter{figure}{0}
\setcounter{equation}{0}
\setcounter{footnote}{0}
\setcounter{table}{0}
\setcounter{section}{0}



\title{$H_D$-Quantum Vertex
  Algebras}

\runningheads{M. J. Bergvelt}{$H_D$-Quantum Vertex
  Algebras}

\begin{start}


\author{M. J. Bergvelt}{1},

\address{Department of Mathematics, University of Illinois\\
  Urbana-Champaign, Illinois 61801, USA}{1}


\begin{Abstract}
We discuss a class of quantum vertex algebras where not only the
commutativity of vertex algebra is broken by a braiding map $S^{(\tau)}$, but also
the translation covariance is broken by a translation map $S^{(\gamma)}$. The new
class of quantum vertex operators satisfy a Braided Jacobi Identity
containing both the braiding and the translation maps.
\end{Abstract}
\end{start}


\section{What is a Vertex Algebra?}
According to Borcherds \cite{MR1653021}, \cite{MR1865087} a \textbf{vertex
algebra} can be thought of as a commutative, associative, unital
\emph{singular} algebra with infinitesimal translation
symmetry. Usually this made precise as follows. We start with a
\emph{State Space} $V$ (a vector space) with a distinguished element 1
(called the \emph{vacuum}, playing the r\^ ole of unit in $V$). $V$
has an action of the Hopf algebra $H_D=\mathbb{C}[D]$, where $D$ is
the infinitesimal translation operator. Then for each $a\in V$ we have
a (singular) operator $Y(a,z)$ of left multiplication by $a$:
\[
Y(a,z)\colon V\to V((z)).
\]
We have the following axioms:
\begin{itemize}
\item \textbf{Vacuum: } $Y(1,z)=1_V$, and $Y(a,z)1=e^{ZD}a$.
\item \textbf{Translation Covariance: } $[D,Y(a,z)]=Y(Da,z)=\partial_z Y(a,z)$.
\item \textbf{Commutativity: }For all $a, b\in V$ there is an $N$ such that
\[
(z_1-z_2)^N[Y(a,z_1),Y(b,z_2)]=0.
\]
\end{itemize}
In this formulation in the multiplication $Y(a,z)b$ the left factor
$a$ is treated differently than the right factor $b$; for instance $a$
has the variable $z$ associated to it, but $b$ has no variable
attached. We find it useful to reformulate the theory so that both
factors are treated symmetrically. We attach to the two factors $a$
and $b$ two variables, $z_1,z_2$, and define a singular multiplication
    \begin{align*}
X_{z_1,z_2}\colon V\ten V&\to V[[z_1,z_2]][(z_1-z_2)\inv],\\
a\ten b&\to X_{z_1,z_2}(a\ten b)=e^{z_2D}Y(a,z_1-z_2)b.
    \end{align*}
Then one derives from the axioms for $Y$ the following properties of $X_{z_1,z_2}$:
\begin{itemize}
  \item \textbf{Vacuum: }$X_{z_1,z_2}(a\ten 1)=e^{z_1 D}a$, $X_{z_1,z_2}(1\ten
    a)=e^{z_2 D}a$. 
  \item \textbf{Translation Covariance: }$e^{\gamma D}X_{z_1,z_2}(a\ten b)=
    X_{z_1+\gamma,z_2+\gamma}(a\ten b)$.
  \item \textbf{Commutativity: }$X_{z_1,z_2}(a\ten b)=X_{z_2,z_1}(b\ten a)$.
\end{itemize}
Relation of $X_{z_1,z_2}$ to $Y(a,z)$ is given by expansions.
Let $i_{z_1;z_2}$ be the {expansion} of a rational function in $z_1-z_2$
in the region $\abs{z_1}>\abs{z_2}$. 
Then
\begin{equation}
\begin{aligned}
i_{z_1;z_2}X_{z_1,z_2}(a\ten b)&= Y(a,z_1)Y(b,z_2)1,\\
i_{z_2;z_1}X_{z_1,z_2}(a\ten b)&= Y(b,z_2)Y(a, z_1)1,\\
i_{z_2;z_3}X_{z_2+z_3,z_2}(a\ten b)&= Y(Y(a,z_3)b,z_2)1.
\end{aligned}\label{eq:expVertexOps}
\end{equation}
From this one easily derives the  axioms for vertex
algebras in terms of $Y$, starting with the properties of  $X_{z_1,z_2}$.

\section{Quantum Vertex Algebras via Deformation.}

As a motivation for our definition of a quantum vertex algebra
consider  a commutative algebra  $M$ , with a multiplication
\[
m\colon M\ten M\to M,
\]
so that we have in particular $m(a\ten b)=m(b\ten a)$. Now a
quantization of $(M,m)$ could be defined by introducing a formal
variable $t$, and a deformed multiplication
\[
m_t\colon M_t\ten M_t\to M_t,\quad M_t=M[[t]],
\]
where $m=m_t\mod t$. In general $m_t$ will not longer be
commutative. We could require that there is a \emph{Braiding Map}:
\[
S\colon M_t\ten M_t\to M_t\ten M_t,
\]
such that $m_t \circ S(a\ten b)= m_t(b\ten a)$. So the braiding
describes the failing of $m_t$ to be commutative. Now if the
commutative algebra $M$ has some symmetry, say via a group $G$ acting
on $M$ (so that $m(ga\ten gb)=gm(a\ten b)$), then we can expect that
after quantisation the deformed multiplication is no longer
$G$-symmetric. Instead, we can require that there is for each $g\in G$
a map
\[
S^g\colon M_t\ten M_t \to M_t\ten M_t,
\]
such that 
\[
g m_t\circ S^g(a\ten b)=m_t(ga\ten gb).
\]
So a deformation of a commutative algebra with symmetry $(M,m,G)$
would be a quintuple $(M_t,m_t,S, S^g, G)$, satisfying a complicated
system of axioms we don't want to write down here, see
\cite{AngBergv}. 

Now a vertex algebra is a singular analog of a commutative algebra
with as symmetry the Hopf algebra $H_D$ of infinitesimal
translations. So if we quantize we can expect that the commutativity
and translation covariance are no longer exact, and that we need extra
structures to describe the broken symmetries. 

Introduce a quantum variable $t$ and deformed singular multiplication
\begin{align*}
X_{z_1,z_2}\colon V^{\ten2}&\to V[[z_1,z_2]][z_1\inv,(z_1-z_2)\inv][[t]].
    \end{align*}
The commutativity and translation covariance are supposed to be not
longer exact. The extra structure we need is
\begin{itemize}
\item a \emph{Braiding} map: 
\[
S^{(\tau)}_{z_1,z_2}\colon V^{\ten2}\to
  V^{\ten2}[z_1^{\pm 1},z_2^{\pm 1},(z_1-z_2)\inv][[t]].
\]
\item a \emph{Translation} map: 
\[
S^{(\gamma)}_{z_1,z_2}\colon 
  V^{\ten2}\to V^{\ten2}[z_1^{\pm
    1},z_2,(z_1-z_2)\inv,(z_1+\gamma)^{\pm1}, z_2+\gamma][[t]].
\]
\end{itemize}
 The deformed multiplication and the braiding and translation maps are
supposed to have the following properties:
\begin{itemize}
\item Braided Commutativity:
\[
X_{z_1,z_2} (a\ten b)=X_{z_2,z_1}\circ S^{(\tau)}_{z_2,z_1}(b\ten a).
\]
\item (Broken) Translation covariance:
\[
e^{\gamma D}X_{z_1,z_2}\circ S^{(\gamma)}_{z_1,z_2}(a\ten
b)=X_{z_1+\gamma,z_2+\gamma}(a\ten b).
\]
\item Plus a bunch of other axioms (Yang-Baxter, hexagon, \dots), see \cite{AngBergv})
\end{itemize}

This defines a \textbf{$H_D$-quantum vertex algebra}.
\section{Vertex Operators and the Braided Jacobi Identity.}

We define a 1-variable vertex operator as usual:
\[
Y(a,z)b=X_{z,0}(a\ten b),
\]
i.e., by evaluating the second variable of $X$ at zero. (Note that in
general in an $H_D$-quantum vertex algebra we can not evaluate the
\emph{first} variable at 0.) Then the relation between $Y$ and $X$ is
given by a variant of \eqref{eq:expVertexOps}, where we need to insert
braiding and translation matrices in appropriate places. More
generally, one shows that there exists 
\[
X_{z_1,z_2,z_3}\colon V^{\ten 3}\to
V[[z_k]][z_i\inv,(z_i-z_j)\inv][[t]],\quad {1\le i<j \le 3}, i\le k\le 3
\]
such that
\[
i_{z_1;z_2,z_3}X_{z_1,z_2,z_3}=X_{z_1,0}(1\ten
X_{z_2,z_3}).
\]
Then we have the following expansions: if $A=a\ten b\ten c$, then
\begin{equation}
\begin{aligned}
    i_{z_1;z_2}X_{z_1,z_2,0}(A)&=Y(a,z_1)Y(b,z_2)c, \\
    i_{z_2;z_1}X_{z_1,z_2,0}(A)&=Y_{z_2}(1\otimes Y_{z_1})
    i_{z_2;z_1}S^{(\tau), 12}_{z_2,z_1}
    (b\otimes a\otimes c),\\
    i_{z_2;z_3}X_{z_2+z_3,z_2,0}(A)&=Y_{z_2}(Y_{z_3}\otimes 1)
    i_{z_2;z_3}S^{(z_2),12}_{z_3,0} (a\otimes b\otimes
    c).
  \end{aligned}\label{eq:ExpansQVertex}
  \end{equation}
Here we write $Y_z(a\ten b)$ for $Y(a,z)b$. From these expansions one
derives the \emph{Braided Jacobi Identity} for $H_D$-quantum vertex
algebras:
    \begin{multline*}
      i_{z_1;z_2}\delta(z_1-z_2,z_3)a(z_1)b(z_2)c-
      i_{z_2;z_1}\delta(z_1-z_2,z_3)Y_{z_2}(1\otimes Y_{z_1})
      S^{(\tau),12}_{z_2,z_1}(b\otimes a\otimes c)\\
      =i_{z_2;z_3}\delta(z_1,z_2+z_3)Y_{z_2}(Y_{z_3}\otimes
      1)S^{(z_2),12}_{z_3,0}(a\otimes b\otimes c).
    \end{multline*}
Here we write $a(z)$ for $Y(a,z)$.

\section{Example: Hall-Littlewood polynomials.}

The main inspiration for our construction came from the theory of
quantum vertex operators for Hall-Littlewood polynomials introduced by
Jing, \cite{MR1112626}. These quantum vertex operators occur naturally
in a $H_D$-quantum vertex algebra $V_{L,t}$ which is a deformation of
the lattice vertex algebra based on the lattice $L=\mathbb{Z}\alpha$,
with pairing $L\ten L\to \mathbb{Z}$, $a\alpha\ten b\alpha\mapsto
ab$. $V_{L,t}$ is in a sense generated by $e^\alpha$ and the braiding
and translation maps in this case are given by
\[
  S^{(\tau)}_{z_1,z_2}(e^{\alpha}\otimes
  e^{\alpha})=-\frac{1-tz_2/z_1}{1-tz_1/z_2}e^\alpha\ten e^\alpha,\quad
  S^{(\gamma)}_{z_1,z_2}(e^{\alpha}\otimes
  e^{\alpha})=\frac{1-tz_2/z_1}{1-t\frac{z_2+\gamma}{z_1+\gamma}}e^\alpha\ten
  e^\alpha,
\]
An effective method to do calculations in $V_{L,t}$ and similar
$H_D$-quantum vertex algebras is given by the
theory of bicharacters, see \cite{Anguelova:thesis}.
\section{Conclusion and Outlook.}

The $H_D$-quantum vertex algebras introduced above are generalizations
of the quantum vertex operators of Etingof-Kazhdan
\cite{MR2002i:17022}. In their theory the vertex operators are
translation covariant, so that the translation maps $S^{(\gamma)}$
are the identity. 

Vertex algebras are commutative algebras with translation covariance
and singularities in the product of vertex operators of the form
$(z_1-z_2)^{-N}$. In our $H_D$-quantum vertex algebras the
commutativity and translation covariance is broken via nontrivial
braiding and translation maps $S^{(\tau)}$ and $S^{(\gamma)}$, but the
singularities are essentially still of the same type
$(z_1-z_2)^{-N}$. Now in examples of quantum vertex operators, see
e.g., \cite{q-alg/9706023}, one sees that in practice the product can
have singularities of the form
\[
\frac{1}{z_1-p^kq^\ell z_2}.
\]
This means that in such quantum vertex algebras one needs to
\emph{extend} the symmetry algebra $H_D=\mathbb{C}[D]$ by adding
(group like) operators $T_p, T_q$ that act like
\[
T_p f(z)=f(pz),\quad T_q f(z)=f(qz).
\]
Note that $H_{p,q}=\mathbb{C}[T_p^{\pm 1},T^{\pm 1}_q,d]$ is non
commutative. 

Now the basic formalism of vertex algebras (and $H_D$-quantum vertex
algebras) is very much based on $H_D$. For instance, the delta
distribution $\delta(z_1,z_2)$, which is ubiquitous in the theory, is
the difference of two expansions of the basic singularity $\frac
1{z_1-z_2}$, 
and expansions are given by the \emph{exponential
  operators} canonically associated to $H_D$.

Replacing the commutative Hopf algebras $H_D$ by the non commutative
$H_{p,q}$ changes the basic framework of vertex algebras drastically:
for instance, one needs a new theory of Dirac delta distributions
adapted to $H_{p,q}$, \cite{AngBergvChiral}.

\section*{Acknowledgments}
The results in this paper are joint work with Iana Anguelova, see
\cite{AngBergv}.



\begin{thebibliography}{Ang06}

\bibitem[ABa]{AngBergv}
I.~I. Anguelova and M.~J. Bergvelt, \emph{{$H_D$}-quantum vertex algebras and
  bicharacters}, preprint, arXiv:math.QA/0706.1528.

\bibitem[ABb]{AngBergvChiral}
\bysame, \emph{{$H_t$}-quantum vertex algebras and deformed chiral algebras},
  in preparation.

\bibitem[Ang06]{Anguelova:thesis}
I.~I. Anguelova, \emph{Bicharacter constructions of quantum vertex algebras},
  Ph.D. thesis, University of Illinois, Urbana-Champaign, 2006.

\bibitem[Bor98]{MR1653021}
Richard~E. Borcherds, \emph{Vertex algebras}, Topological field theory,
  primitive forms and related topics (Kyoto, 1996), Progr. Math., vol. 160,
  Birkh\"auser Boston, Boston, MA, 1998, pp.~35--77. \MR{MR1653021 (99m:17034)}

\bibitem[Bor01]{MR1865087}
\bysame, \emph{Quantum vertex algebras}, Taniguchi Conference on Mathematics
  Nara '98, Adv. Stud. Pure Math., vol.~31, Math. Soc. Japan, Tokyo, 2001,
  pp.~51--74. \MR{MR1865087 (2002k:17054)}

\bibitem[EK00]{MR2002i:17022}
Pavel Etingof and David Kazhdan, \emph{Quantization of {L}ie bialgebras. {V}.
  {Q}uantum vertex operator algebras}, Selecta Math. (N.S.) \textbf{6} (2000),
  no.~1, 105--130. \MR{2002i:17022}

\bibitem[FR97]{q-alg/9706023}
Edward Frenkel and Nikolai Reshetikhin, \emph{{Towards Deformed Chiral
  Algebras}}, Proceedings of the Quantum Group Symposium at the XXIth
  International Colloquium on Group Theoretical Methods in Physics, Goslar
  1996, 1997, arXiv:q-alg/9706023.

\bibitem[Jin91]{MR1112626}
Nai~Huan Jing, \emph{Vertex operators and {H}all-{L}ittlewood symmetric
  functions}, Adv. Math. \textbf{87} (1991), no.~2, 226--248. \MR{MR1112626
  (93c:17039)}

\end{thebibliography}
\def\cprime{$'$}
\providecommand{\bysame}{\leavevmode\hbox to3em{\hrulefill}\thinspace}
\providecommand{\MR}{\relax\ifhmode\unskip\space\fi MR }
\providecommand{\MRhref}[2]{%
  \href{http://www.ams.org/mathscinet-getitem?mr=#1}{#2}
}
\providecommand{\href}[2]{#2}








\end{document}